\documentclass[12pt]{amsart}
\usepackage{amssymb}
\usepackage[mathscr]{eucal}
\usepackage{epsfig}
\usepackage{verbatim}
\usepackage{graphicx}
\usepackage{hyperref}

\usepackage{vmargin}
\setpapersize{USletter}
\setmargrb{1in}{1in}{1in}{1in}

\newcommand{\C}{\ensuremath{\mathbb C}}

\newcommand{\PP}{\ensuremath{\mathbb P}}

\DeclareMathOperator{\rk}{rk}
\DeclareMathOperator{\HB}{HB}

\DeclareMathOperator{\scriptHom}{Hom}

\DeclareMathOperator{\Sing}{Sing}
\DeclareMathOperator{\Zeros}{Zeros}

\theoremstyle{definition}
\newtheorem{defn}{Definition}[section]

\newtheorem{example}[defn]{Example}

\theoremstyle{plain}
\newtheorem{lem}[defn]{Lemma}

\newtheorem*{lem*}{Lemma}
\newtheorem{thm}[defn]{Theorem}
\newtheorem*{thm*}{Theorem}
\newtheorem{cor}[defn]{Corollary}
\newtheorem{prop}[defn]{Proposition}

\theoremstyle{remark}
\newtheorem{rem}[defn]{Remark}

\title[Intersection of curves through points in $\PP^2$]{On the intersection of the curves through a set of points in $\PP^2$}
\author{Zachariah C.~Teitler}
\date{5/1/06}

\address{Department of Mathematics, SLU 10687, Hammond, LA 70402}
\email{zteitler@selu.edu}

\begin{document}

\bibliographystyle{plain}       % Set the bibliography style to AMS
                                % alphabetized. (Can use ``amsalpha'' or
                                % ``abbrv''instead.)

\begin{abstract}
Given a set of points in $\PP^2$,
we consider the common zeros of the set of curves of a given degree
passing through those points.
For general sets of points, these zero sets
have the expected dimension and are smooth.
In fact, given graded Betti numbers,
for any arrangement of points whose ideal
has those graded Betti numbers, general among such arrangements,
the zero sets have the expected dimension and are smooth.
\end{abstract}

\maketitle

\section{Introduction}

There has been a great deal of interest in the linear series
of curves in $\PP^2$ containing a given set of points
(see, for example,
\cite{gambier}, 
\cite{MR760027},
\cite{MR846019},
\cite{MR1673756},
or \cite{MR0024165}).
In this paper, we consider the intersection of all the curves
of a given degree containing
a given set of points in $\PP^2$.

Let $Z \subset \PP^2$ be an arrangement of points in $\PP^2$ and
$I$ the homogeneous ideal of $Z$.
By ``arrangement'' we mean a finite set of points.
Write $I_d$ for the degree $d$ piece of $I$.
If $d \gg0$, then of course $\Zeros(I_d) = Z$.
We ask: what can one say about $\Zeros(I_d)$ for values of $d$ smaller
than the generating degree of $I$?
For example:
What is the dimension of $\Zeros(I_d)$?
Is it smooth?
The answers to these questions depend partly on the resolution type
of the ideal $I$.
We give answers for arrangements which are general of a given resolution type.

Recall that a finite set $Z$ of points in $\PP^2$
is defined by a Hilbert--Burch matrix,
a matrix whose entries are homogeneous forms on $\PP^2$,
and this matrix determines the minimal free resolution of $I$
(see section~\ref{subsect:res-type}).
Recall also that there are integers $k$, $0 < a_1\leq \dots \leq a_{k+1}$,
$0 < b_1 \leq \dots b_k$
such that the $(i,j)$th entry of the Hilbert--Burch matrix
has degree $b_j-a_i$.
In fact, the $a_i$ are exactly the degrees of the generators
in a minimal generating set of $I$.

Suppose we are given a resolution type as follows.
Let us be given some $(a_1,\dots,a_{k+1};b_1,\dots,b_k)$ such that $b_j > a_i$ for every $i,j$.
Consider the set of arrangements $Z$ defined by Hilbert--Burch matrices
whose entries have degree $b_j-a_i$.
The requirement $b_j > a_i$ means that
for the ideal $I$ of an arrangement $Z$ in this set,
every relation (syzygy) of $I$ has higher degree than every generator of $I$.
For general arrangements $Z$ in this set, we are able to
give answers to the questions above.
Explicitly, we prove the following:
\begin{thm*}
Let us fix $k$, $\{a_i\}$, $\{b_j\}$ as above, such that every $b_j>a_i$.
Consider the set of arrangements defined by Hilbert--Burch matrices
whose $(i,j)$th entries have degree $b_j-a_i$.
Let $Z$ be a general arrangement in this set, and let $I$ be the ideal of $Z$.
Then for $d\geq 0$, $\Zeros(I_d)$ is smooth and has the expected dimension.
Here, expected dimension means the following:
If $d < a_1$ then $\Zeros(I_d) = \PP^2$.
If $a_1 \leq d < a_2$ then $\Zeros(I_d)$ is a curve.
If $a_2 \leq d < a_3$ then $\Zeros(I_d)$ is a finite set.
If $a_3 \leq d$ then $\Zeros(I_d) = Z$.
\end{thm*}
(See Theorem~\ref{thm:main}.)

In particular, for any $n>0$, we give explicit information for
general arrangements of $n$ points, see Corollary~\ref{cor:general-d-envelopes}.

For simplicity, we work over $\C$, but any algebraically closed field of
characteristic zero will do.
The restriction on characteristics comes from the use of Kleiman's
generic smoothness theorem~\cite[III.10.7]{hartshorne},
in the proof of Proposition~\ref{prop:general-transversality}.

\section*{Acknowledgments}

I would like to thank my advisor, Robert Lazarsfeld.
I would also like to thank Mel Hochster for his helpful suggestions
regarding Proposition~\ref{prop:decomposition-of-det-loci}.

\section{Plane arrangements of points}

We introduce terminology for the objects of study,
the intersection of the curves of a given degree through a given set of points.
We also consider families of point arrangements and resolution data.

\subsection{Degree envelopes}

\begin{defn}
Let $Z \subset \PP^n$ be a non-empty closed subscheme with homogeneous ideal $I$.
For $d\geq 0$, we define the \textbf{$d$th degree envelope},
or \textbf{$d$-envelope}, of $Z$ to be
the closed subscheme $Z_d = \Zeros(I_d) \subset \PP^n$
given by the intersection of all the degree $d$
hypersurfaces containing $Z$.
The degree envelopes form a decreasing chain which begins with $\PP^n$
and stabilizes at $Z$.
If $Z_d \neq Z_{d-1}$, we say $d$ is a \textbf{geometric generating degree} of $I$.
\end{defn}

Equivalently, $Z_d$ is the base scheme of the linear series of degree $d$
hypersurfaces containing $Z$.

\begin{example}\label{ex:envelopes}
\begin{enumerate}
\item If $Z$ is a complete intersection of type $(d_1,\dots,d_r)$
with $d_1 < \dots < d_r$, then the geometric generating degrees of $I$
are exactly the $d_i$.
For each $i$, let $H_i$ be a hypersurface of degree $d_i$ such that
$Z = H_1 \cap \dots \cap H_r$.
Then $Z_{d_1} = H_1$, $Z_{d_2} = H_1 \cap H_2$, and so on.
\item Let $Z$ be five general reduced points in $\PP^2$.
Then $Z_2$ is the unique conic containing $Z$, and $Z_3 = Z$.
The geometric generating degrees are $2$ and $3$.
\item Let $Z$ be eight general reduced points in $\PP^2$.
Then there is a pencil of cubics passing through $Z$,
so $Z_3$ consists of the nine basepoints of this pencil.
That is, $Z_3$ is the union of $Z$ with an extra ninth point
(distinct from $Z$ because $Z$ is general).
The geometric generating degrees are $3$ and $4$.
\item\label{ex:3-collinear} Let $Z$ be four reduced points in $\PP^2$ with three collinear, but not all four.
Say the points $P_1$, $P_2$, and $P_3$ lie on the line $L$, and the point $P_4$
lies off of $L$.
Then $Z_2 = L \cup P_4$ and $Z_3 = Z$.
In this case a degree envelope has components of different dimensions.
The geometric generating degrees are $2$ and $3$.
\item\label{ex:eleven}
Let $C$ be a smooth plane cubic and let $Z$ be eleven general reduced points on $C$.
Then $Z_3 = C$.
There is a unique point $P \in C$ such that $Z \cup P$ is the complete intersection
of $C$ with a quartic curve, and $Z_4 = Z \cup P$, twelve points ($P$ is distinct from
all the points of $Z$ by generality).
Finally, $Z_5 = Z$.
In this case, $I$ has three geometric generating degrees, $3$, $4$, and $5$.
\item
Let $Z$ be a set of $18$ points in $\PP^2$ in general position.
Then the ideal $I(Z)$ is minimally generated by three forms of degree $5$
and one form of degree $6$, but only $5$ is a geometric generating degree of $I(Z)$.
That is, $6$ is a degree of a generator of $I(Z)$, but not a geometric generating degree.
\end{enumerate}
\end{example}

The following lemma will clarify the relationship between
the geometric generating degrees of $I$ and the usual
degrees of (algebraic) generators of $I$.

\begin{lem}\label{lem:ggds}
Let $I \subset S = \C[x_0,\dots,x_n]$ be a saturated homogeneous ideal.
Say $I = (H_1,\dots,H_s)$, with $H_i$ homogeneous, $\deg H_i = d_i$,
and $d_1 \leq \dots \leq d_s$. Then:
\begin{enumerate}
\item every geometric generating degree of $I$ is one of the integers $d_i$.
\item $d_1$ is a geometric generating degree of $I$.\qed
\end{enumerate}
\end{lem}

\begin{rem}
We regard $\PP^n$ and $Z$ itself as trivial degree envelopes of $Z$.
(We have $Z_d = \PP^n$ for $d < d_1$, in the notation of Lemma~\ref{lem:ggds},
and $Z_d = Z$ for $d \gg 0$.)
So $Z$ has no non-trivial degree envelopes if and only if $I$ has
only one geometric generating degree.
\end{rem}

These degree envelopes arise naturally in the following situation.
Let us consider an arrangement of lines through the origin of $\C^3$.
Let $A\subset \C^3$ be the union of these lines and let $I$
be the homogeneous ideal of $A$.
If we blow up the origin, then the total transform of the ideal $I$
may have embedded components supported in the exceptional divisor
of the blowup.
The exceptional divisor is a $\PP^2$ on which the strict transforms of the lines
in $A$ mark out an arrangement of points.
It is shown in the companion paper~\cite{zct:mla}
that the non-trivial degree envelopes (in $\PP^2$) of this point arrangement
are the supports of embedded components of the total transform of $I$.
The geometric generating degrees of $I$ determine
the structure of these embedded components.

This situation arose in the process of computing the multiplier
ideals of such an ideal of an arrangement of lines in $\C^3$,
as explained in~\cite{zct:mla}.
Corollary~\ref{cor:general-d-envelopes} is used in that paper
to discuss general arrangements of lines.

\subsection{Partition of $(\PP^2)^n$ by graded Betti numbers}\label{subsect:res-type}

The Hilbert--Burch theorem gives a useful description of the defining ideal of
a Cohen--Macaulay subvariety of codimension $2$ in a smooth projective variety.
(See, for example, \cite{burch:hb-thm} or \cite[Section 20.4]{eisenbud:comm-alg}.)
A configuration of finitely many points in $\PP^2$ is the first example of such
a subvariety.
We state the theorem only in this special case.

\begin{thm}[Hilbert--Burch]\label{thm:hb}
Let $Z \subset \PP^2$ be a finite set
(a zero-dimensional reduced closed subscheme)
with saturated homogeneous ideal $I \subset S = \C[x,y,z]$.
Then there is an integer $k>0$ and integers $0 < a_1 \leq \cdots \leq a_{k+1}$,
$0 < b_1 \leq \cdots \leq b_k$ such that
the minimal graded free resolution of $I$ has the form,
\begin{equation}\label{eqn:hb-resolution-algebraic}
  0 \to \bigoplus_{j=1}^k S(-b_j)
    \overset{A}{\longrightarrow} \bigoplus_{i=1}^{k+1} S(-a_i)
    \to I \to 0 ,
\end{equation}
where $A$ is a $(k+1)\times k$ matrix of homogeneous forms.
The ideal $I$ is generated by the determinants of the $k\times k$ minors of $A$.
\end{thm}
\begin{proof}
See, for example, \cite[Theorem 4.3]{eisenbud:syzygies}.
\end{proof}
The $a_i$ and $b_j$ are the \textbf{resolution data} of $I$.
One can verify $\sum a_i = \sum b_j$.
The resolution data is equivalent to the \textbf{graded Betti numbers}
of $I$~\cite{eisenbud:syzygies}.
To be precise, the graded Betti numbers give in degree $d$ the number of times
that $d$ occurs on the lists $\{a_i\}$ and $\{b_j\}$.

\begin{defn}\label{defn:resolution-data}
\textbf{Resolution data} is a pair of lists $(\{a_i\}_{i=1}^{k+1},\{b_j\}_{j=1}^{k})$
with $0 < a_1 \leq \dots \leq a_{k+1}$, $0 < b_1 \leq \dots \leq b_k$,
and $\sum a_i = \sum b_j$.

We say resolution data $R=(\{a_i\},\{b_j\})$ is \textbf{positive}
if $a_{k+1} < b_1$ (so that $a_i < b_j$ for every $i$ and $j$).
\end{defn}

\begin{rem}\label{rem:number-of-points}
Let $Z$ be an arrangement of $n$ points
with resolution data $R=(\{a_i\},\{b_j\})$, not necessarily positive.
One can show that $n=(\sum b_j^2 - \sum a_i^2)/2$, for example by
computing the dimensions of global sections of large twists
of the short exact sequence~\eqref{eqn:hb-resolution-algebraic}.
See also~\cite[Exercise 3.15]{eisenbud:syzygies}.
\end{rem}

The collection of all arrangements of $n$ distinct points on $\PP^2$ corresponds
naturally to $(\PP^2)^n - \Delta$, the open complement of the diagonals in $(\PP^2)^n$,
up to choosing an ordering for the $n$ points.
This open set is partitioned by resolution data
(equivalently, by graded Betti numbers)
into pieces that are constructible sets in the Zariski topology \cite{MR1219271}.

The main goal of this paper is to prove the following.

\begin{thm}\label{thm:main}
Let $(\{a_i\}_{i=1}^{k+1},\{b_j\}_{j=1}^{k})$ be positive resolution data.
Let $T \subset (\PP^2)^n$ be the locus of arrangements
with this resolution data, where $n=(\sum b_j^2 - \sum a_i^2)/2$
as in~\ref{rem:number-of-points}.
Then $T$ is irreducible.
Let $Z \in T$ be a general arrangement.
Let $I$ be the ideal of $Z$.
Then the set of geometric generating degrees of $I$ and
the degree closures of $Z$ are as follows.

If $k=1$, then the geometric generating degrees of $I$ are $\{a_1,a_2\}$.
In particular, $Z_d = \PP^2$ for $d < a_1$ and $Z_d = Z$ for $d \geq a_2$.
For $a_1 \leq d < a_2$, $Z_d = Z_{a_1}$ is smooth with codimension $1$.

If $k\geq 2$, then the geometric generating degrees of $I$ are $\{a_1,a_2,a_3\}$.
In particular, $Z_d = \PP^2$ for $d < a_1$ and $Z_d = Z$ for $d \geq a_3$.
For $a_1 \leq d < a_2$, $Z_d=Z_{a_1}$ is smooth with codimension $1$
and for $a_2 \leq d < a_3$, $Z_d=Z_{a_2}$ is smooth with codimension $2$
(that is, a set of reduced points).
\end{thm}

\begin{rem}
The case $k=1$ in Theorem~\ref{thm:main} corresponds to complete intersections.
\end{rem}
%\item Part~\ref{ex:eleven} of Example~\ref{ex:envelopes} corresponds
%to the resolution data
%\[
%  a_1 = 3, a_2 = 4, a_3 = 5 ; \quad b_1 = b_2 = 6 .
%\]
%\item Consider the resolution data
%\[
%  a_1 = a_2 = 4, a_3 = 6; \quad b_1 = b_2 = 7 .
%\]
%A general arrangement $Z$ with this resolution data consists of
%of $15$ of the points of a (transversal) complete intersection
%of type $(4,4)$.
%Then $Z_4=Z_5$, the sixteen points of the complete intersection,
%and $Z_6=Z$.

\begin{rem}
In Example~\ref{ex:envelopes}(\ref{ex:3-collinear}), the resolution data
is $(2,2,3;3,4)$, hence not positive.
Note that the $2$-envelope consists of a line plus a point,
so this fails the codimension part of the conclusion of the theorem.

In general it is not known what happens when the points have non-positive resolution data.
\end{rem}

As a special case, so to speak, we get explicit information for
general sets of $n$ points in $\PP^2$, meaning all arrangements
corresponding to points in some fixed open subset of $(\PP^2)^n$.

\begin{cor}\label{cor:general-d-envelopes}
Let $n>1$.
Let $Z$ be a set of $n$ general points in $\PP^2$.
Let $d$ and $r$ be specified by $\binom{d+1}{2} \leq n = \binom{d+2}{2} - r$ with $r>0$,
so that $d$ is the lowest degree of a curve passing through $Z$,
and $r$ is the number of independent curves of degree $d$ passing through $Z$.
Let $I$ be the ideal of $Z$.
\begin{enumerate}
\item If $r=1$, the geometric generating degrees of $I$ are $\{d,d+1\}$
and the $d$-envelope $Z_d$ is a smooth curve of degree $d$.
\item If $r=2$ and $d>2$, the geometric generating degrees of $I$ are $\{d,d+1\}$,
and $Z_d$ is a set of $d^2$ distinct,
reduced points in $\PP^2$, a complete intersection of type $(d,d)$,
containing $Z$ together with $d^2-n = \binom{d-1}{2}$ extra points.
\item If $r=2$ and $d=2$ (so $n=4$),
then $2$ is the only
geometric generating degree of $I$.
\item If $r\geq 3$, then $d$ is the only
geometric generating degree of $I$.
\end{enumerate}
\end{cor}

\begin{proof}
It suffices to note that the partition of $(\PP^2)^n$ by graded Betti numbers
includes a dense piece, corresponding to certain resolution data
given in~\cite{geramita-gregory-roberts}.
We repeat this ``generic'' resolution data here.
Let $r$ and $d$ be defined as in the statement of the theorem.
Then, with notation as in Theorem~\ref{thm:hb},
the ``generic'' values of $k$, $\{a_i\}$, and $\{b_j\}$ are as follows.
\begin{itemize}
\item If $2r \geq d+2$ then $k=d+1-r$, $a_1=\dots=a_{k+1} = d$,
$b_1=\dots=b_{2r-d-2} = d+1$, and $b_{2r-d-1}=\dots=b_k=d+2$.
\item If $2r \leq d+2$ then $k=d+1-r$, $a_1=\dots=a_r=d$,
$a_{r+1}=\dots=a_{k+1}=d+1$, and $b_1=\dots=b_k=d+2$.
\end{itemize}
A general arrangement of $n$ points has this resolution data,
and we apply Theorem~\ref{thm:main}.
If $r \geq 3$, then $a_1=a_2=a_3=d$, so $d$ is the only geometric generating
degree of the ideal $I$ of the arrangement.
The other cases $r=1,2$ are similar.
\end{proof}

To prove Theorem~\ref{thm:main}, we interpret an arrangement $Z$
and its Hilbert--Burch matrix in terms of a vector bundle and
apply general transversality results.

\section{Arrangements via  vector bundles}

In this section we reinterpret point arrangements in $\PP^2$
and their degree envelopes
in terms of sections of a vector bundle.

\subsection{The Hilbert--Burch vector bundle}

Given resolution data $R=(\{a_i\},\{b_j\})$
as in Definition~\ref{defn:resolution-data},
we define the \textbf{Hilbert--Burch vector bundle}
\[
  E(R) = \scriptHom \left( \bigoplus \mathscr{O}_{\PP^2}(-b_j) ,
      \bigoplus \mathscr{O}_{\PP^2}(-a_i) \right) .
\]
Note, $E(R)$ is ample if and only if $R$ is positive.
We define the \textbf{vector space of Hilbert--Burch matrices of type $R$}
to be $\HB(R) = H^0(\PP^2,E(R))$,
the vector space of $(k+1)\times k$ matrices whose $(i,j)$th
entry is a homogeneous form of degree $b_j-a_i$ for each $i$, $j$.
For $A \in \HB(R)$, let $I(A)$ be the ideal generated by the determinants
of the $k\times k$ minors of $A$
and let $Z(A)\subset\PP^2$ be the subscheme cut out by $I(A)$.

\begin{thm}\label{thm:general-hb-mx}
Let $R = (a_1,\dots,a_{k+1};b_1,\dots,b_k)$
be positive resolution data as in Definition~\ref{defn:resolution-data}.
Let $A \in \HB(R)$ be general.
\begin{enumerate}
\item\label{thm:already-known} $Z(A)$ is an arrangement of
distinct, reduced points, with resolution data $R$.
The number of points is $n = (\sum b_j^2 - \sum a_i^2)/2$.
\item For each $d \geq 0$, the $d$-envelope $Z(A)_d$
is smooth with codimension determined as follows.
Let $r(d) = \#\{a_i \leq d\}$.
Explicitly, $r(d)$ is defined by $a_{r(d)} \leq d < a_{r(d)+1}$,
with $r(d)=0$ for $d<a_1$, and $r(d)=k+1$ for $d\geq a_{k+1}$.
Then $Z(A)_d$ has codimension $r(d)$ if $r(d)\leq 2$.

Furthermore, if $k=1$, then $Z(A)_d = Z(A)$ if and only if $r(d) = 2$;
if $k \geq 2$, then $Z(A)_d=Z(A)$ if and only if $r(d) > 2$.
\end{enumerate}
\end{thm}

\begin{rem}
Part~\ref{thm:already-known} is already well-known.
\end{rem}

This easily implies Theorem~\ref{thm:main}.

\begin{proof}[Proof of Theorem~\ref{thm:main}]
Let $R = (\{a_i\},\{b_j\})$ and $T$ be as in the statement of
Theorem~\ref{thm:main}.
By the first part of Theorem~\ref{thm:general-hb-mx},
the map $A \mapsto Z(A)$
is a rational map
$\HB(R) \dashrightarrow T$.
It is surjective, by the Hilbert--Burch theorem~\ref{thm:hb}.
Since $\HB(R)$ is irreducible, it follows that $T$ is irreducible.

For general $Z \in T$, there is a (general) $A \in \HB(R)$
such that $Z = Z(A)$.
Then the claims of Theorem~\ref{thm:main} regarding $Z$
follow immediately from Theorem~\ref{thm:general-hb-mx} applied to $A$.
\end{proof}

To prove Theorem~\ref{thm:general-hb-mx}, we interpret
the degree envelopes $Z(A)_d$ as loci where $A$, as a section of the
Hilbert--Burch bundle, meets certain cones.
The rest of this section is devoted to developing these tools,
and then at the end we prove the theorem.

\subsection{Decomposition of determinantal loci}

Let $X$ be a generic $(k+1)\times k$ matrix of variables whose entries $x_{ij}$
are independent variables.
For $1 \leq i \leq k+1$, let $F_i$ be the determinant of the $k\times k$
minor of $X$ obtained by deleting the $i$th row.

\begin{defn}\label{defn:det-loci}
Let $S = \C[x_{1,1},\dots,x_{k+1,k}]$.
Let $M = M_{(k+1)\times k}(\C)$, the vector space of $(k+1)\times k$ matrices
with constant entries.
The entries $x_{ij}$ of $X$ give coordinates on $M$.
For $1 \leq r \leq k+1$, we define certain ideals and determinantal loci in $M$,
as follows.
\begin{enumerate}
\item Let $I_r \subset S$ be the ideal $(F_1,\dots,F_r)$.
\item Let $J_r \subset S$ be the ideal generated by the determinants of
the maximal minors of the $(k+1-r)\times k$ matrix consisting of the last
$k+1-r$ rows of $X$ (all but the first $r$ rows).
In particular, we set $J_{k+1} = (1)$.
\item Let $L_r \subset M$ be the subscheme cut out by $I_r$.
\item Let $N_r \subset M$ be the subscheme cut out by $J_r$.
\end{enumerate}
\end{defn}

By a theorem of
Eagon and Hochster~\cite{hochster:generic-determinantal-ideals-are-radical},
$I_{k+1}$ is prime, as are all the $J_r$.
So $L_{k+1}$ is irreducible, and so are all the $N_r$.
We have the following very useful decomposition of the determinantal loci $L_r$:
\begin{prop}\label{prop:decomposition-of-det-loci}
For $1 \leq r \leq k+1$, $L_r$ is reduced,
and $L_r = L_{k+1} \cup N_r$ as schemes.
\end{prop}
I am grateful to M.~Hochster for suggesting to me the proof of this statement.
It follows from
the statement on ideals that $I_r = I_{k+1} \cap J_r$,
which we prove shortly.

\begin{example}
Say $k=2$, so
\[ X = \begin{pmatrix} a & b \\ c & d \\ e & f \end{pmatrix} . \]
We have
\[ F_1 = cf - de , \quad F_2 = af - be , \quad F_3 = ad - bc , \]
so
\[
  I_1 = (cf - de) , \quad I_2 = (cf - de , af - be) ,
  \quad I_3 = (cf - de , af-be , ad-bc) ,
\]
and
\[
  J_1 = (cf-de) , \quad J_2 = (e,f) , \quad J_3 = (1) .
\]
Obviously $I_1 = I_3 \cap J_1$ and $I_3 = I_3 \cap J_3$,
but we also have, less obviously, $I_2 = I_3 \cap J_2$.
\end{example}

\begin{lem}
Let $R$ be a ring, $I \subset R$ be an ideal, and $e \notin I$.
Assume:
$P = I + (e)$ is radical,
$Q = (I:e) = \{\, x \mid xe \in I \,\}$ is prime,
and $e^2 \notin I$ (equivalently, $e \notin Q$).
Then $I = P \cap Q$.
\end{lem}
\begin{proof}
First we show $I$ is radical.
Suppose $x^n \in I$ for $n\geq 2$.
Then $x^n \in P$, so $x \in P$.
Therefore $x = i + ae$ for some $i \in I$.
Since $x^n \in I$ and $i \in I$, we get $(ae)^n \in I$;
in particular, $a^n e^{n-1} \in Q$.
Since $Q$ is prime and $e \notin Q$, $a \in Q$.
Thus $ae \in I$, so $x \in I$.

Now, suppose $y \in P \cap Q$.
We may write $y = i + ae$, with $i \in I$.
Then $y^2 = iy + aye$, where $iy \in I$ and $ye\in I$ because $y \in Q$.
Therefore $y^2 \in I$.
Since $I$ is radical, $y \in I$.
\end{proof}

\begin{proof}[Proof of Proposition~\ref{prop:decomposition-of-det-loci}]
We go by downward induction on $r$, starting from $r=k+1$.
Since $J_{k+1} = (1)$, the unit ideal, the initial case is trivial.

For $r \leq k$, $I_{r+1} = I_r + (F_{r+1})$; this ideal is radical by induction.
We claim that $J_r = (I_r : F_{r+1})$ and $F_{r+1} \notin J_r$.
From these claims and the previous lemma it follows
that $I_r = I_{r+1} \cap J_r$, and in particular that $I_r$ is radical.

For the second claim, note that
\[ F_{r+1} \notin (x_{r+1,1}, \dots, x_{r+1,k}) \supset J_r . \]

For the first claim,
if $G F_{r+1} \in I_r \subset J_r$
then, since $J_r$ is prime and $F_{r+1} \notin J_r$,
we have $G \in J_r$.
This shows $(I_r : F_{r+1}) \subset J_r$.

We have to show $J_r F_{r+1} \subset I_r$.
We claim that for any generator $P$ of $J_r$ given as a maximal minor
of the last $k+1-r$ rows of $X$, we have $P F_{r+1} \in I_r$.
We may reorder the columns of $X$ so that the minor whose determinant gives $P$
is given by the first $k+1-r$ columns of $X$.
Take the transpose of these columns and write it in block form as $(UV)$,
where $U$ is the first $r$ columns and $V$ is the square matrix of size $k+1-r$
whose determinant is $P$.
Let $w=(F_1,-F_2,\dots,(-1)^i F_i,\dots)$,
and write it also in block form as $w=(w_1,w_2)$ where $w_1$ has size $r$
and $w_2$ has size $k+1-r$.
Then by Cramer's rule,
\[  \begin{pmatrix} U & V \end{pmatrix} \begin{pmatrix} w_1 \\ w_2 \end{pmatrix} = 0 , \]
so $V w_2 = - U w_1$.
Multiplying on the left by the adjoint matrix $V^{*}$ of $V$ (the transpose
of the matrix of cofactors) and applying again Cramer's rule,
\[ P w_2 = \det(V) w_2 = -V^{*}U w_1 . \]
In particular, $F_{r+1}$ is the first entry of $w_2$,
so $P F_{r+1}$ is some combination of the entries of $w_1$,
namely $F_1,\dots,F_r$.
This shows $P F_{r+1} \in I_r$.
Therefore $J_r F_{r+1} \subset I_r$, and so $(I_r : F_{r+1}) = J_r$.

Applying the previous lemma, we see that $I_r = I_{r+1} \cap J_r$, and by induction,
\[ I_r = I_{r+1} \cap J_r = I_{r+2} \cap J_{r+1} \cap J_r = I_{r+3} \cap J_{r+2} \cap J_{r+1} \cap J_r = \dots \]
Since $J_r \subset J_{r+1} \subset \dots$, we see that, as claimed,
$ I_r = I_{k+1} \cap J_r $.
\end{proof}

We will take advantage of the following useful facts about $L_{k+1}$ and the $N_r$.
\begin{prop}\label{fact:codims-of-det-loci}
\begin{enumerate}
\item $L_{k+1}$ has codimension $2$ in $M$.
\item The singular locus $\Sing L_{k+1}$ has codimension $6$ in $M$.
\item Each $N_r$ has codimension $r$ in $M$.
\item Each $\Sing N_r$ has codimension $2(r+1)$ in $M$.
\item $L_{k+1} \subset N_1 = L_1$, but $L_{k+1} \not\subset N_r$ for any $r>1$.
\item $N_{k+1} = \emptyset \subset L_{k+1}$, but $N_r \not\subset L_{k+1}$
for any $r < k+1$.
\item For $1 < r < k+1$, $L_{k+1} \cap N_r$ has codimension at least $3$ in $M$.
\end{enumerate}
\end{prop}
\begin{proof}
We use the well-known formula that in the space of $m \times n$ matrices,
the variety of matrices with rank at most $c$
has codimension equal to $(m-c)(n-c)$ (see, for example,~\cite[Prop.~12.2]{harris:intro}),
and singular locus equal to the variety of matrices with rank at most $c-1$
(see, for example,~\cite[Example~14.16]{harris:intro}).
We apply this to prove the first four parts as follows.

For~(1), $L_{k+1}$ is the variety of matrices with rank at most $k-1$
in the space of $(k+1)\times k$ matrices.
For~(2), $\Sing L_{k+1}$ is the variety of matrices with rank at most $k-2$,
in the same space.

Now, write $M = M_1 \times M_2$, where $M_1$ is the affine space with coordinates
given by the entries of the first $r$ rows of $X$, and $M_2$ is the affine space
with coordinates given by the last $k+1-r$ rows of $X$.
Let $N'_r \subset M_2$ be the locus defined by the vanishing of all the maximal
minors of the last $k+1-r$ rows of $X$.
Then
\[ N_r = M_1 \times N'_r . \]
Since $N'_r$ is the variety of matrices of rank at most $k-r$
in the space of $(k+1-r) \times k$ matrices, $N'_r$ has codimension $r$ in $M_2$.
Therefore $N_r$ has codimension $r$ in $M$, proving~(3).
Furthermore,
\[ \Sing N_r = M_1 \times \Sing N'_r , \]
where $\Sing N'_r \subset M_2$ has codimension $2(r+1)$.
This proves~(4).

For~(5) and~(6), the inclusions $L_{k+1} \subset N_1$ and $N_{k+1} \subset L_{k+1}$
are clear.
To see the noninclusions, consider the following $(k+1)\times k$ matrices,
given in block form:
\[
  A_r = \begin{pmatrix} 0 & 0 \\ I_{k+1-r} & 0 \end{pmatrix},
  \quad
  B = \begin{pmatrix} I_k \\ 0 \end{pmatrix}
\]
where $I_{k+1-r}$ and $I_k$ are the identity matrices of the indicated sizes.
Then for $r>1$, $A_r \in L_{k+1}$ but $A_r \notin N_r$.
For $r < k+1$, $B \in N_r$ but $B \notin L_{k+1}$.

Finally, for~(7), for $1 < r < k+1$, $L_{k+1} \cap N_r$ is strictly contained in $L_{k+1}$,
which is irreducible and of codimension $2$.
\end{proof}

\subsection{Cones in $E(R)$ and degree envelopes}\label{subsect:cones}

Given positive resolution data $R$ and a Hilbert--Burch matrix $A \in \HB(R)$,
recall that $I(A)$ is the ideal of determinants of $k\times k$ minors
of $A$.
Each of these is obtained by omitting a row of $A$.
For $1 \leq i \leq k+1$, let $F_i(A)$ be the determinant of the $k\times k$
minor of $A$ obtained by omitting the $i$th row.
Note that $\deg F_i(A) = a_i$.
Then $I(A) = (F_1(A),\dots,F_{k+1}(A))$,
and for $d \geq 0$, the $d$-envelope $Z(A)_d$
is defined by the forms $F_i(A)$ such that $\deg F_i(A) = a_i \leq d$.

The matrix $A$ is also a section of the Hilbert--Burch bundle $E(R)$,
and we take advantage of this to give an alternative approach for
$Z(A)$ and its degree envelopes.
The idea is to define cones in the total space of $E(R)$ analogous
to the $L_r \subset M$ considered in the previous section,
and then recover $Z(A)$ and the $Z(A)_d$ as the loci in $\PP^2$
where $A$ meets these cones.

We denote by $\mathbf{E}(R)$ the total space of the
vector bundle $E(R)$.
Let $\pi:\mathbf{E}(R) \to \PP^2$ be the projection map.
There is a tautological map of bundles on $\mathbf{E}(R)$,
\[
  \pi^{*} \bigoplus_{j=1}^k \mathscr{O}_{\PP^2}(-b_j)
    \to \pi^{*} \bigoplus_{i=1}^{k+1} \mathscr{O}_{\PP^2}(-a_i) .
\]
Abusing notation, we denote this tautological map by $X$,
and for each $i,j$, we denote by $x_{ij}$ the induced map
\[ x_{ij} : \pi^{*} \mathscr{O}_{\PP^2}(-b_j) \to \pi^{*} \mathscr{O}_{\PP^2}(-a_i) . \]
The $x_{ij}$ are global coordinates on $\mathbf{E}(R)$.
Suppose over an affine open subset $U \subset \PP^2$ one trivializes
each of the line bundles $\mathscr{O}_{\PP^2}(-b_j)$, $\mathscr{O}_{\PP^2}(-a_i)$.
We get a trivialization of $E(R)|_U$, hence coordinates
on $\mathbf{E}(R)|_U = \pi^{-1}(U)$.
These coordinates are the $x_{ij}$ (together with coordinates on $U$).
In particular, the $x_{ij}$ restrict to coordinates on each fiber of $\mathbf{E}(R)$.

We define cones in $\mathbf{E}(R)$ by vanishing of determinants of minors
of $X = (x_{ij})$,
just as in the previous section.
As before, for $1\leq r \leq k+1$, let $F_r$ be the determinant
of the minor of $X$ obtained by omitting the $r$th row.
The vanishing-locus $\{F_r=0\}\subset \mathbf{E}(R)$ is the rank-dropping locus
of the vector bundle map given by removing the $r$th row of $X$:
\[
  \pi^{*}\bigoplus_{j=1}^k \mathscr{O}_{\PP^2}(-b_j)
    \to \pi^{*}\bigoplus_{\substack{1\leq i \leq k \\ i\neq r}} \mathscr{O}_{\PP^2}(-a_i) .
\]
For $1 \leq r \leq k+1$, let $L_r(R) \subset \mathbf{E}(R)$
be defined by $F_1=\dots=F_r=0$, the scheme-theoretic intersection
of the rank-dropping loci.

Similarly, let $N_r(R) \subset \mathbf{E}(R)$ be defined by the vanishing
of the maximal minors of the last $k+1-r$ rows of $X$.
Equivalently, $N_r$ is the rank-dropping locus of the map of vector bundles,
\[
  \pi^{*}\bigoplus_{j=1}^k \mathscr{O}_{\PP^2}(-b_j)
    \to \pi^{*}\bigoplus_{i=r+1}^{k+1} \mathscr{O}_{\PP^2}(-a_i) .
\]

Now, over an affine open $U \subset \PP^2$, trivializing each $\mathscr{O}_{\PP^2}(-b_j)$,
$\mathscr{O}_{\PP^2}(-a_i)$, the resulting trivialization of $E(R)|_U$ gives an isomorphism
\[ \mathbf{E}(R)|_U \longrightarrow M \times U , \]
which takes
\begin{align*}
  L_r(R)|_U \longrightarrow L_r \times U , \\
  N_r(R)|_U \longrightarrow N_r \times U
\end{align*}
This leads to the following ``global'' analogue of
Propositions~\ref{prop:decomposition-of-det-loci} and~\ref{fact:codims-of-det-loci}:
\begin{prop}\label{prop:global-det-loci}
Let $R$ be positive resolution data.
For $1 \leq r \leq k+1$, $L_r(R)$ is reduced and $L_r(R) = L_{k+1}(R) \cup N_r(R)$.
Also, $L_{k+1}(R)$ is irreducible and reduced, and each $N_r(R)$ is irreducible and reduced.
We have the following facts:
\begin{enumerate}
\item $L_{k+1}(R)$ has codimension $2$ in $\mathbf{E}(R)$.
\item $\Sing L_{k+1}(R)$ has codimension $6$ in $\mathbf{E}(R)$.
\item Each $N_r(R)$ has codimension $r$ in $\mathbf{E}(R)$.
\item Each $\Sing N_r(R)$ has codimension $2(r+1)$ in $\mathbf{E}(R)$.
\item $L_{k+1}(R) \subset N_1(R) = L_1(R)$,
but $L_{k+1}(R) \not\subset N_r(R)$ for any $r>1$.
\item $N_{k+1}(R) = \emptyset \subset L_{k+1}(R)$,
but $N_r(R) \not\subset L_{k+1}(R)$
for any $r < k+1$.
\item For $1 < r < k+1$, $L_{k+1}(R) \cap N_r(R)$ has codimension
at least $3$ in $\mathbf{E}(R)$.\qed
\end{enumerate}
\end{prop}

We have defined the cones we are interested in.
Now we want to show how to use them to get point arrangements in $\PP^2$
and degree envelopes.

For positive resolution data $R$ and $A \in \HB(R)$,
the arrangement $Z(A)$ and its degree envelopes $Z(A)_d$ are defined by
the vanishing of the forms $F_i(A)$.
The idea is to see these $F_i(A)$ as pullbacks of the equations $F_i$
on $\mathbf{E}(R)$, and then we will see that the $Z(A)$ and $Z(A)_d$ are
the loci in $\PP^2$ where $A$, as a section of $\mathbf{E}(R)$,
intersects the cones $L_r(R)$.

Let $s_A:\PP^2 \to \mathbf{E}(R)$ be the section associated to $A$.
For $1 \leq i \leq k+1$ and $1 \leq j \leq k$,
one has the following two maps of line bundles:
\begin{gather*}
  A_{ij}:\mathscr{O}_{\PP^2}(-b_j) \to \mathscr{O}_{\PP^2}(-a_i) , \\
  x_{ij}: \pi^{*} \mathscr{O}_{\PP^2}(-b_j)
                \to \pi^{*} \mathscr{O}_{\PP^2}(-a_i).
\end{gather*}
Evidently $A_{ij} = s_A^{*} x_{ij}$.
This implies, for $1 \leq i \leq k+1$, $F_i(A) = s_A^{*} F_i$.
We obtain the following:
\begin{prop}
Let $R=(\{a_i\},\{b_j\})$ be positive resolution data and $A \in \HB(R)$.
Let $s_A:\PP^2 \to \mathbf{E}(R)$ be the map
corresponding to $A\in H^0(\PP^2,E(R))$.
Then $Z(A) = s_A^{-1}(L_{k+1}(R))$, the locus in $\PP^2$ where
$A$ meets $L_{k+1}(R)$.

For $d \geq 1$, the $d$-envelope $Z(A)_d$
is the scheme-theoretic preimage $s_A^{-1}(L_r(R))$,
the locus where $A$ meets $L_r(R)$,
where $r=r(d)=\#\{a_i \leq d\}$ (or, $r$ is defined by $a_r \leq d < a_{r+1}$---the
same function $r(d)$ as in Theorem~\ref{thm:general-hb-mx}).\qed
\end{prop}

\subsection{General transversality for sections of a vector bundle}

We recall the following well-known statement:
\begin{lem}
Let $E$ be a globally generated vector bundle of rank $e$ on a projective variety $X$.
Let the total space of $E$ be denoted $\mathbf{E}$.
Let $L \subset \mathbf{E}$ be a closed subset
with $\dim L < \rk E$.
Then a general section $s \in H^0(X,E)$ does not meet $L$.\qed
\end{lem}
This is proved by a dimension count.
It can be generalized to give the following proposition,
reminiscent of the proof of Bertini's
theorem in characteristic zero via Kleiman's generic smoothness theorem
as presented in~\cite[III.10.9]{hartshorne}.
It belongs to the folklore, but for lack of a reference we give a statement
and proof.

\begin{prop}\label{prop:general-transversality}
Let $E$ be an ample and globally generated vector bundle of rank $e$
on a smooth complex projective variety $X$.
Let the total space of $E$ be denoted $\mathbf{E}$.
Let $L \subset \mathbf{E}$ be an irreducible reduced closed subset
with $\dim L \geq e$ and $\dim \Sing L < e$.
Then for a general global section $s$ of $E$,
the locus $s^{-1}(L) \subset X$ where $s$ meets $L$ is
nonempty, reduced, smooth, and with codimension in $X$
equal to the codimension of $L$ in $\mathbf{E}$.
\end{prop}
\begin{proof}
By Theorem~12.1(c) of~\cite{fulton:big}, for every section $s$ of $E$,
$s^{-1}(L)$ is a positive cycle on $X$, so in particular nonempty.

Let $U \subset H^0(X,E)$ be the open subset of sections not meeting $\Sing L$.
Consider
\[ \tilde{L} = \{ \, (s,x) \in U \times X \mid s(x) \in L - \Sing(L) \, \} . \]
This is a nonempty open subset of the set
\begin{equation} \label{eqn:total-space-affine-bundle}
  \{ \, (s,x) \in H^0(X,E) \times X \mid s(x) \in L - \Sing(L) \, \}  .
\end{equation}
Note that 
\begin{equation}  \label{eqn:restricted-bundle}
  \{ \, (s,x) \in H^0(X,E) \times X \mid s(x) \in L - \Sing(L) \, \}  \longrightarrow L-\Sing(L)
\end{equation}
is an affine bundle.
Indeed, it is the restriction to $L-\Sing L$ of the affine bundle
$H^0(X,E) \times X \longrightarrow \mathbf{E}$ given by $(s,x) \mapsto s(x)$.
The restricted bundle~(\ref{eqn:restricted-bundle}) has smooth base $L-\Sing L$;
therefore its total space~(\ref{eqn:total-space-affine-bundle})
is reduced and smooth.
Hence the open subset $\tilde{L}$ is reduced and smooth.

The projection map $\tilde{L} \to U$ is surjective 
because every section of $E$ meets $L$.
By Kleiman's generic smoothness theorem~\cite[Theorem III.10.7]{hartshorne},
there is an open dense subset $W \subset U$ over which
the fibers of this projection map are nonempty, reduced, smooth,
and all of the same codimension,
namely the codimension of $L$ in $\mathbf{E}$.
Finally, the fiber over $s \in W$ is isomorphic to $s^{-1}(L) \subset X$.
\end{proof}

\subsection{Proof of Theorem~\ref{thm:general-hb-mx}}

We now use the tools we have just developed to prove
Theorem~\ref{thm:general-hb-mx},
in turn implying
Theorem~\ref{thm:main}
and Corollary~\ref{cor:general-d-envelopes}.

\begin{proof}[Proof of Theorem~\ref{thm:general-hb-mx}]
Let $R=(\{a_i\}_{i=1}^{k+1},\{b_j\}_{j=1}^{k})$
be the positive resolution data given in the
hypothesis of the statement of the theorem.
Let $E(R)$ be the Hilbert--Burch vector bundle
as defined above, with total space $\mathbf{E}(R)$,
and for $1 \leq r \leq k+1$
let $L_r(R), N_r(R) \subset \mathbf{E}(R)$
be the cones defined in section~\ref{subsect:cones}.
The positivity of $R$ means $E(R)$ is a direct sum of ample
line bundles on $\PP^2$, hence ample.

We saw that for every $A \in \HB(R)$,
with corresponding section $s_A : \mathbb{P}^2 \to \mathbf{E}(R)$,
the subscheme $Z(A) \subset \PP^2$
is the locus where $s_A$ meets $L_{k+1}(R)$.
Recall that by Proposition~\ref{prop:global-det-loci},
$L_{k+1}(R)$ has codimension $2$ in $\mathbf{E}(R)$
and $\Sing L_{k+1}(R)$ has codimension $6$ in $\mathbf{E}(R)$.
Therefore, for general sections $A \in \HB(R)$,
 Proposition~\ref{prop:general-transversality}
shows that $Z(A)$ is nonempty, reduced, and smooth,
with codimension $2$ in $\PP^2$.
One checks easily that the Hilbert--Burch short exact sequence
as in Theorem~\ref{thm:hb} is a resolution of the ideal $I(A)$,
so $Z(A)$ has the resolution data $R$, as claimed.
The number of points is $n=(\sum b_j^2 - \sum a_i^2)/2$ of $Z(A)$
by an argument as in Remark~\ref{rem:number-of-points}.
This proves the first part of the theorem, which was nevertheless
known previously.

Now let $d \geq 0$.
Recall our earlier notation, that for $A \in \HB(R)$ and $1 \leq r \leq k+1$,
$F_r(A)$ is the homogeneous form of degree $a_r$ given by the determinant
of the $k\times k$ minor of $A$ obtained by omitting the $r$th row.
Then the $d$-envelope $Z(A)_d$ is defined by the vanishing of those forms
$F_r(A)$ such that $\deg F_r(A) = a_r \leq d$.
Since $a_1 \leq \cdots \leq a_{k+1}$,
we see that $Z(A)_d$ is defined by $F_1(A)=\dots=F_r(A)=0$,
where $r=r(d)=\#\{a_i \leq d\}$ is given, as in the statement of the theorem,
by $a_r \leq d < a_{r+1}$,
with $r=0$ for $d<a_1$ and $r=k+1$ for $d \geq a_{k+1}$.

First of all, if $r=0$, then $Z(A)_d = \PP^2$ is clear.

Suppose $r \geq 1$.
We saw in section~\ref{subsect:cones} that $Z(A)_d$
is the locus where the corresponding section $s_A : \PP^2 \to \mathbf{E}(R)$ meets $L_r(R)$.
We now apply Proposition~\ref{prop:global-det-loci}, as follows.

If $r=1$ (equivalently, $a_1\leq d < a_2$),
then $L_1(R) = N_1(R)$, which is irreducible and of codimension $1$
in $\mathbf{E}(R)$, with singularities $\Sing(N_1(R))$ of codimension $4$ in $\mathbf{E}(R)$.
Then by Proposition~\ref{prop:general-transversality},
for general $A \in \HB(R)$, $Z(A)_d$
is smooth and reduced, with codimension $1$.
(Note that $Z(A)_d$ is defined by the single equation $F_1(A)$.)

If $r=2$ (equivalently, $a_2 \leq d < a_3$),
then $L_2(R) = L_{k+1}(R) \cup N_2(R)$.
In this case $L_{k+1}(R)$ and $N_2(R)$ both have codimension $2$ in $\mathbf{E}(R)$.
Since $L_{k+1}(R) \cap N_2(R)$ has codimension at least $3$,
by Proposition~\ref{prop:global-det-loci}, we see that
for general $A \in \HB(R)$, the section $s_A$ does not meet $L_{k+1}(R) \cap N_2(R)$.
Therefore for such $A$, $Z(A)_d$ is the disjoint union
\[ Z(A)_d = s_A^{-1}(L_{k+1}(R)) \cup s_A^{-1}(N_2(R)) = Z(A) \cup s_A^{-1}(N_2(R)) . \]
Since $N_2(R)$ has codimension $2$ in $\mathbf{E}(R)$ and
$\Sing N_2(R)$ has codimension $6$, Proposition~\ref{prop:general-transversality}
shows $s_A^{-1}(N_2(R))$ is nonempty, smooth, reduced, and of codimension $2$ in $\PP^2$.
Therefore $Z(A)_d$ is a reduced set of points, strictly larger than $Z(A)$.

If $r \geq 3$ (equivalently, $a_3 \leq d$),
then $L_r(R) = L_{k+1}(R) \cup N_r(R)$.
Since $N_r(R)$ has codimension $r > \dim \PP^2$ in $\mathbf{E}(R)$,
general sections $A \in \HB(R)$ do not meet $N_r(R)$.
Therefore the $d$-envelope $Z(A)_d$, the locus where $s_A$ meets $L_r(R)$,
is just the locus where $s_A$ meets $L_{k+1}(R)$.
This is $Z(A)$.
Therefore for $d \geq a_3$, $Z(A)_d = Z(A)$.

This completes the proof of the theorem.
\end{proof}

%\subsection{Further questions}
\begin{rem}
It is natural to consider similar questions in higher dimension.
One expects similar results for general arrangements of points in $\mathbb{P}^s$:
that the geometric generating degrees and degree envelopes are determined
by the number $n$ of points.

One may also consider more special arrangements of points in $\PP^s$.
For example, Gorenstein point arrangements in $\PP^3$ are defined by the Pfaffians
of a skew-symmetric matrix (see~\cite{MR0453723}),
and certain point arrangements in $\PP^s$ are defined by the
minors of a $k\times (k+s-1)$ matrix.
\end{rem}

\bibliography{../biblio}

\begin{thebibliography}{10}

\bibitem{MR0453723}
David~A. Buchsbaum and David Eisenbud.
\newblock Algebra structures for finite free resolutions, and some structure
  theorems for ideals of codimension {$3$}.
\newblock {\em Amer. J. Math.}, 99(3):447--485, 1977.

\bibitem{burch:hb-thm}
Lindsay Burch.
\newblock On ideals of finite homological dimension in local rings.
\newblock {\em Proc. Cambridge Philos. Soc.}, 64:941--948, 1968.

\bibitem{MR1219271}
S.~Diaz and A.~Geramita.
\newblock Points in projective space in very uniform position.
\newblock {\em Rend. Sem. Mat. Univ. Politec. Torino}, 49(2):267--280 (1993),
  1991.

\bibitem{eisenbud:comm-alg}
David Eisenbud.
\newblock {\em Commutative algebra}, volume 150 of {\em Graduate Texts in
  Mathematics}.
\newblock Springer-Verlag, New York, 1995.
\newblock With a view toward algebraic geometry.

\bibitem{eisenbud:syzygies}
David Eisenbud.
\newblock {\em The geometry of syzygies}, volume 229 of {\em Graduate Texts in
  Mathematics}.
\newblock Springer-Verlag, New York, 2005.
\newblock A second course in commutative algebra and algebraic geometry.

\bibitem{fulton:big}
William Fulton.
\newblock {\em Intersection theory}, volume~2 of {\em Ergebnisse der Mathematik
  und ihrer Grenzgebiete. 3. Folge. A Series of Modern Surveys in Mathematics
  [Results in Mathematics and Related Areas. 3rd Series. A Series of Modern
  Surveys in Mathematics]}.
\newblock Springer-Verlag, Berlin, second edition, 1998.

\bibitem{gambier}
Bertrand Gambier.
\newblock Syst{\`e}me lin{\'e}aire de courbes alg{\'e}briques de degr{\'e}
  donn{\'e} admettant un groupe donn{\'e} de points bases.
\newblock {\em Annale Scientifiques de l'{\'E}cole Normale Sup{\'e}rieure},
  S{\'e}r.~3(41):147--264, 1924.

\bibitem{geramita-gregory-roberts}
A.~V. Geramita, D.~Gregory, and L.~Roberts.
\newblock Monomial ideals and points in projective space.
\newblock {\em J. Pure Appl. Algebra}, 40(1):33--62, 1986.

\bibitem{MR760027}
A.~V. Geramita and P.~Maroscia.
\newblock The ideal of forms vanishing at a finite set of points in {${\bf
  P}\sp{n}$}.
\newblock {\em J. Algebra}, 90(2):528--555, 1984.

\bibitem{MR846019}
Brian Harbourne.
\newblock The geometry of rational surfaces and {H}ilbert functions of points
  in the plane.
\newblock In {\em Proceedings of the 1984 Vancouver conference in algebraic
  geometry}, volume~6 of {\em CMS Conf. Proc.}, pages 95--111, Providence, RI,
  1986. Amer. Math. Soc.

\bibitem{harris:intro}
Joe Harris.
\newblock {\em Algebraic geometry}, volume 133 of {\em Graduate Texts in
  Mathematics}.
\newblock Springer-Verlag, New York, 1995.
\newblock A first course, Corrected reprint of the 1992 original.

\bibitem{hartshorne}
Robin Hartshorne.
\newblock {\em Algebraic geometry}.
\newblock Springer-Verlag, New York, 1977.
\newblock Graduate Texts in Mathematics, No. 52.

\bibitem{hochster:generic-determinantal-ideals-are-radical}
M.~Hochster and John~A. Eagon.
\newblock Cohen-{M}acaulay rings, invariant theory, and the generic perfection
  of determinantal loci.
\newblock {\em Amer. J. Math.}, 93:1020--1058, 1971.

\bibitem{MR1673756}
Rick Miranda.
\newblock Linear systems of plane curves.
\newblock {\em Notices Amer. Math. Soc.}, 46(2):192--201, 1999.

\bibitem{MR0024165}
Louis Nollet.
\newblock Recherches sur les syst\`emes lin\'eaires de courbes alg\'ebriques
  planes.
\newblock {\em M\'em. Soc. Roy. Sci. Li\'ege (4)}, 7:469--554, 1947.

\bibitem{zct:mla}
Zachariah~C. Teitler.
\newblock Multiplier ideals of general line arrangements in $\mathbb{C}^3$.
\newblock arXiv:math.AG/0508308, 2005.

\end{thebibliography}

\end{document}